\documentclass[a4paper]{amsart}

\usepackage[utf8]{inputenc}
\usepackage{amsmath}
\usepackage{amssymb}
\usepackage{amsthm}
\usepackage{tikz}
\usepackage{graphicx}
\usepackage{float}
\usepackage{svg}
\usepackage{verbatim}
\usepackage{url}
\usepackage{mathtools}
\usepackage{microtype}
\usepackage{thm-restate}
\usepackage{enumitem}
\usepackage[hidelinks]{hyperref}
\usepackage[nameinlink]{cleveref}
\usepackage{caption}
\let\cref\Cref
\Crefname{enumi}{}{}

\theoremstyle{plain}
\newtheorem{theorem}{Theorem}[section]
\theoremstyle{definition}
\newtheorem{definition}[theorem]{Definition}
\theoremstyle{plain}

\theoremstyle{definition}
\newtheorem{example}[theorem]{Example}

\theoremstyle{plain}

\newtheorem{corollary}[theorem]{Corollary}
\theoremstyle{remark}
\newtheorem{remark}[theorem]{Remark}

\newcommand{\qua}{\hskip 0.4em \ignorespaces}
\def\arxiv#1{\relax\ifhmode\unskip\qua\fi
\href{http://arxiv.org/abs/#1}%
{\tt arXiv:\penalty -100\unskip#1}}

\def\MR#1{\relax\ifhmode\unskip\qua\fi
\href{https://mathscinet.ams.org/mathscinet-getitem?mr=#1}{\tt MR#1}}
\def\ZB#1{\relax\ifhmode\unskip\qua\fi
\href{https://zbmath.org/?q=an:#1}{\tt Zbl\:#1}}
\def\xox#1{\csname xx#1\endcsname}

\renewenvironment{thebibliography}[1]{
  \begin{oldthebibliography}{#1}\small
    \setlength{\itemsep}{.5ex}
    \setlength{\parskip}{0em}
}
{
  \end{oldthebibliography}
}

\let\stdthebibliography\thebibliography
\let\stdendthebibliography\endthebibliography
\renewenvironment*{thebibliography}[1]{%
  \stdthebibliography{MM20.}}
  {\stdendthebibliography}

\newcommand{\myemail}[1]{\href{mailto:#1}{#1}}

\newcommand{\hfhat}{\widehat{HF}}

\title{On binding sums of contact manifolds}
\address{ETH Z\"urich, R\"amistrasse 101, 8092 Z\"urich, Switzerland\bigskip}
\author{Miguel Orbegozo Rodriguez}
\email{\myemail{miguel.orbegozorodriguez@math.ethz.ch}}
\urladdr{\url{https://people.math.ethz.ch/~morbegozo/}}
\begin{document}
%\subjclass[2020]{57K33, 57K18 (57K20)}

\begin{abstract}
In this short note, we give examples of binding sums of contact 3--manifolds that do not preserve properties such as tightness or symplectic fillability. We also prove vanishing of the Heegaard Floer contact invariant for an infinite family of binding sums where the summands are Stein fillable. This recovers a result of Wendl and Latschev--Wendl~\cite{Wendl, Latschev-Wendl}. Along the way, we rectify a subtle computational error in~\cite{JuhaszKang} concerning the spectral order of a neighbourhood of a Giroux torsion domain.
\end{abstract}

\maketitle

\section{Introduction}
There are several operations that may be performed to obtain new manifolds from existing ones, and the main questions regarding these operations are whether the resulting manifold inherits properties from the original ones. For what is one of the most common operations, the connected sum, this is indeed the case; and properties of the manifolds originating in this manner are well understood.

In this note we concern ourselves with an operation on contact $3$--manifolds which is somewhat similar to the connected sum, which we call the \emph{binding sum}, and show that in general it does not preserve properties of contact structures. We will thus work throughout in the contact category.

While the connected sum of two manifolds consists of removing two balls from them and gluing along the $S^2$ boundary, the binding sum consists of removing solid tori, and then gluing along the torus boundary. For our purposes, the solid tori will be neighbourhoods of binding components of supporting open books; in particular, the construction will depend on the choice of open book. See \Cref{sec:def} for a precise definition and \Cref{sec:examples} for examples.

Several properties of connected sums can be determined from properties of the original manifolds. For instance, the Heegaard Floer homology of a connected sum is the tensor product of the homologies of its summands, and moreover this also translates to the contact setting, where the Heegaard Floer invariant behaves in the same way. Indeed, for contact manifolds $(M_i, \xi_i)$, we get $c(\xi_1 \# \xi_2) = c(\xi_1) \otimes c(\xi_2)$.
Even more, the connect sum inherits the weakest of the properties of the two manifolds, e.g.~the connect sum of a Stein fillable manifold with a strongly fillable one is strongly fillable. 

For the binding sum, this does not happen in general. Indeed, it is possible to do a binding sum of Stein fillable manifolds such that the resulting manifold is supported by an open book that is overtwisted, which in particular means it is not fillable in any sense (not even tight). We provide several examples of this.

However, in these examples the summands have multiple binding components but the sum is performed along a single one. A natural question is whether similar pathological phenomena can occur when we sum along all binding components. Once again, the answer is yes. We show that for an infinite family of sums such that the summands are Stein fillable (in particular their contact invariant is nonzero) the result of performing a binding sum on them gives a contact manifold with vanishing contact invariant. This is our main result.

\begin{theorem} \label{thm:main}
    Let $\Sigma_g$ be the closed surface of genus $g$, and let $\Sigma_{g,2}$ be the connect sum of $ \Sigma_g $ with an annulus $(S^1\times [0,1])$. Then, for any $g > 0$ we get the following.

    \begin{itemize}
        \item [a)] The binding sum of the Stein fillable open books $(\Sigma_{g,2}, \textrm{Id}) $ and $ (S^1 \times [0,1],\textrm{Id})$, where the sum is performed on both boundary components, is a contact manifold whose contact invariant vanishes.
        \item [b)] The binding sum of the open books $(\Sigma_{g,2}, \varphi)$ and $ (S^1 \times [0,1],\textrm{Id})$, where the sum is performed on both boundary components; and $\varphi$ is a diffeomorphism supported in $\Sigma_g$, is a contact manifold whose contact invariant vanishes.
    \end{itemize}
    
\end{theorem}

 Part $a)$ of \Cref{thm:main} also follows from work of Wendl~\cite{Wendl} and Latschev--Wendl~\cite{Latschev-Wendl}. This result, however, uses $J$--holomorphic techniques, and we want to understand this explicitly from the open book perspective. More precisely, we use Klukas' description of an open book of a binding sum \cite{Klukas} to compute an explicit chain in Heegaard Floer homology that kills the contact class. Moreover, this allows us to prove part $b)$ of \Cref{thm:main}, which, due to the diffeomorphism of $\Sigma_g$, does not fit into Wendl's construction. 
 
 Note once again that the monodromy $\varphi$ can be chosen to support a Stein fillable contact structure, thus providing more examples of Stein fillable open books that give rise to contact manifolds with vanishing contact class when the binding sum is performed.

Moreover, to prove \Cref{thm:main} we calculate vanishing of the contact class for the partial open book of a neighbourhood of a Giroux torsion domain, and in doing so we correct a subtle computational error in \cite{JuhaszKang}. In that paper, Juh\'{a}sz and Kang then use the chain that kills the contact invariant to obtain a bound for the \emph{spectral order} (as defined in \cite{SpectralOrder}) of a manifold with Giroux torsion. Denoting by $\mathbf{o}(M,\xi)$ the spectral order of the contact manifold $(M, \xi)$, they obtain the following.

\begin{corollary} [\cite{JuhaszKang}] \label{cor:spectralbound}
    If a contact manifold $(M, \xi)$ with convex boundary has Giroux $2\pi$--torsion, then $\mathbf{o}(M,\xi) \leq 2$. 
\end{corollary}

With our correction, we are able to prove \Cref{cor:spectralbound}. Moreover, although it is conjectured in \cite{SpectralOrder} that the spectral order should be $1$, we remark that the spectral order bound that we get is the same as the one computed in \cite{JuhaszKang}.

\subsection*{Structure of the paper}
In \Cref{sec:def} we formally define the binding sum, and outline Klukas' construction of an open book of the resulting manifold in terms of the original open books. In \Cref{sec:examples} we give several examples to showcase that the binding sum does not, in general, preserve properties of contact manifolds. In \Cref{sec:nonsym}, we prove \Cref{thm:main}. Finally, in \Cref{sec:spectralorder} we explain the correction of the computation in \cite{JuhaszKang} and prove \Cref{cor:spectralbound}.

\subsection*{Acknowledgements}
Most of the mathematical work presented here was completed during my PhD, and thus I would like to thank my supervisor Andy Wand for support. I would also like to thank Peter Feller for valuable conversations. Finally, I would like to gratefully acknowledge financial support by the SNSF Grant 181199 and by the DFG Grant 513007277.

\section{Preliminaries on binding sums of contact manifolds}

For $i = 1,2$ let $(M_i, \xi_i)$ be closed contact 3-manifolds, and $(B_i, \pi _i)$ open book decompositions supporting them. Let $K_i$ be a component of the binding $B_i$. Remove a standard tubular neighbourhood of these knots, which we can identify as the normal bundle $\nu K_i$, and glue the resulting boundaries via a map that preserves the fibres. Since the contact structure in a neighbourhood of a binding component is standard (more generally, it is standard in a neighbourhood of any transverse knot, see for example \cite{Klukas}), we can glue the contact structures and the resulting manifold inherits a contact structure from the contact structures of the original manifolds.

\begin{definition}
This operation is called a \emph{binding sum} and will be denoted by $M_1 \boxplus_{K_1,K_2} M_2$.
\end{definition}

\begin{remark}
Observe that we can extend this definition to summing along several binding components and not just one. Moreover, to avoid complicating the notation, we will drop the subscripts $K_i$ when it is clear which binding components are involved in the sum.
\end{remark}

\begin{figure}[htp]
    \centering
    \includegraphics[height=4cm]{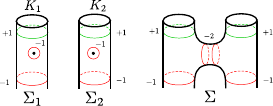}
    \caption[Open book for a binding sum.]{On the left, neighbourhoods of the binding components together with their respective navels. The black dots represent the intersections of the knots $K'_1, K'_2$ with the pages. On the right, the new open book near the binding components being summed. Away from this neighbourhood the monodromy is the composition of the original monodromies.}
    \label{fig:BindingOB}
\end{figure}

We can see from the construction that performing a binding sum interferes with the fibrations of the open book decompositions, so we do not immediately get an open book of the manifold $M_1 \boxplus M_2$. Nevertheless, by a result of Klukas \cite{Klukas} we can obtain an open book decomposition of the summed manifold from the (abstract) open books of the original manifolds. Indeed, given a binding component $K$, let $K'$ be a transverse knot that intersects each page of the open book exactly once near the boundary. Note that since the monodromy is the identity in a neighbourhood of the boundary, this is indeed a knot. Now change the monodromy of the open book near the boundary as follows.
Add a negative Dehn twist along a curve that is the boundary of a neighbourhood of the intersection of $K'$ with the page, a boundary parallel positive twist, and a negative twist along a curve that is boundary parallel but on the other side of the intersection of $K'$ with the page. Call the composition of these Dehn twists $f_K$. Clearly $f_K$ is isotopic to the identity, and thus for any open book $(\Sigma, \varphi)$, the open book $(\Sigma, \varphi \circ f_K)$ is equivalent to it, since $\varphi$ and $\varphi \circ f_K$ represent the same mapping class. In the language of Klukas, we say that $K$ \emph{admits a navel}, and clearly every binding component of an open book admits a navel. The knot $K'$ is then called the \emph{core of the navel}. We can see this on the left hand side of \Cref{fig:BindingOB}. 

\begin{remark}
    Klukas' setup is slightly more general than ours, since he considers knots with arbitrary framing, and we will restrict ourselves to knots with $0$ framing.
\end{remark}

Klukas then shows that a binding component $K$ is transversely isotopic to the core of its navel. Therefore, for open books $(\Sigma_1, \varphi_1)$, $(\Sigma_2, \varphi_2)$ with binding components $K_1, K_2$,  instead of removing a neighbourhood of the binding components $K_1, K_2$ and gluing around the resulting boundary torus, we can remove a neighbourhood of the core of their navels $K_1', K_2'$ and glue around the resulting boundary. Now, on each page this results in removing a neighbourhood of the intersection point of $K_i'$ with the page, and gluing along the boundary circles. This means that the fibrations of the open book decompositions are preserved, so we do get a new open book. The new page $\Sigma$ is the connect sum of the original pages, and the new monodromy is given by the composition of (the natural extensions of) $\varphi_1 \circ f_{K_1}$ and $\varphi_2 \circ f_{K_2}$. Note that $f_{K_1}$ and $f_{K_2}$ are no longer isotopic to the identity after performing the connect sum, see the right hand side of \Cref{fig:BindingOB}.\medskip

\begin{remark}
In the case where we deal with abstract open books (as opposed to ambient ones), we will denote the binding sum as $(\Sigma_1, \varphi_1) \boxplus (\Sigma_2,\varphi_2)$. Again when it is clear which binding components are involved in the sum we will drop them from the notation.
\end{remark}

 \label{sec:def}

\section{Examples} \label{sec:examples}
Now that we have an open book description of the result of a binding sum, we can use it to investigate whether properties of the original manifolds are reflected in the sum. In this section we provide several examples to illustrate that this is not the case even at a basic level.

To begin with, note that for the connect sum --as long as the summands are connected--it does not matter which points we choose to remove neighbourhoods of. However, if we take different supporting open books, we remove different tori, and the result is not necessarily the same contact manifold, and even the underlying smooth manifold may not be the same. This is illustrated in \Cref{ex:S3example} below.

\begin{example} \label{ex:S3example}
    Summing two open books $(D^2,\textrm{Id})$ (representing the standard Stein fillable $S^3)$ gives the manifold determined by the open book $(S^1 \times [0,1],\textrm{Id})$, which is the standard Stein fillable $S^1 \times S^2$. 
    We can see this on the left--hand side of \Cref{fig:S3BindingSum}, where the bottom twists vanish because they are along contractible curves, and the two negative twists cancel with the positive boundary twists.
    
    However, if we stabilize both open books to get Hopf bands and sum those (along one boundary component), we get an open book that is not right-veering and thus overtwisted by \cite{HKM}. We can see this on the right--hand side of \Cref{fig:S3BindingSum}, where the positive boundary twists of the original open books get cancelled with the negative twists coming from the binding sum.
    Moreover, we can destabilize the two boundary twists to end with $(S^1 \times [0,1],\tau^{-2})$, which is an open book for $\mathbb{RP}^3$. Thus, not only do we get an overtwisted contact structure--as opposed to a tight one before-- but we now have a different underlying manifold.

\begin{figure}[htp]
    \centering
    \includegraphics[width = 0.9 \linewidth]{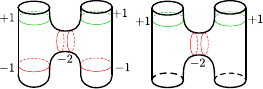}
    \caption[Binding sums depend on the open book.]{On the left, the sum $(D^2, \textrm{Id}) \boxplus (D^2, \textrm{Id}) $ gives $(S^1 \times [0,1],\textrm{Id})$. On the right, $(S^1 \times [0,1], \tau) \boxplus (S^1 \times [0,1], \tau)$ gives an open book that is not right-veering.}
    \label{fig:S3BindingSum}
\end{figure}
\end{example}

Note that overtwistedness of the open book on the right--hand side of \Cref{fig:S3BindingSum} follows from existence of an arc that veers to the left. This is a more general phenomenon. Indeed, if one of the open books has more than one binding component, and we do not perform the binding sum on all boundary components we may obtain an open book that is not right-veering. We can see this in the following \Cref{ex:LVsum}.

\begin{example} \label{ex:LVsum}
If we sum two open books $(\Sigma_1,\textrm{Id})$ and $(\Sigma_2,\varphi)$, where $\Sigma_1$ is any surface which has a boundary component not used in the sum, we get a left-veering arc by connecting the non summed boundary component with a summed one, due to the negative twist at the bottom, see \Cref{fig:OTBindingSum}. Note that $(\Sigma_1, \textrm{Id})$ is Stein fillable (it is the standard Stein fillable contact structure on $\#^k(S^1 \times S^2)$), and the construction does not depend on $\varphi$, so we can choose it to support a Stein fillable contact structure.

\begin{figure}[htp]
    \centering
    \includegraphics[height=4.5cm]{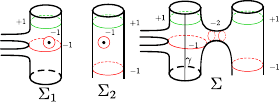}
    \caption[Left-veering arc in a binding sum.]{On the left the two open books $(\Sigma_1,\textrm{Id})$ and $(\Sigma_2,\varphi)$ being summed, together with their corresponding navels. On the right, the open book for the binding sum and a left-veering arc $\gamma$.}
    \label{fig:OTBindingSum}
\end{figure}

\end{example}

Up until now, the examples provided give overtwisted manifolds from Stein fillable ones because we create arcs that veer to the left. This, however, is not the only way to achieve overtwistedness. For instance, stabilize the sum from \Cref{ex:LVsum} enough to make it right-veering (this can always be done by \cite{HKM}). Stabilizations do not change the contact manifold, so we still have an overtwisted contact structure. Then, we observe that, if none of the components being summed was a disc (so the left--veering arc happens in a boundary component that was not used in the sum), the stabilizations are performed away from the part of the open book where we have performed the binding sum and, in particular, these two operations commute.

What is more, it turns out that we do not need to stabilize to get overtwisted, right--veering open books, as the following example shows. 

\begin{example}
Consider the open books $(S^1 \times [0,1],\tau^n)$ and $(S^1 \times [0,1], \tau^2)$, where $n \geq 2$. In particular, both of these open books are non--destabilizable and support Stein fillable contact structures. The result of summing them along one boundary component is a right-veering, non-destabilizable, overtwisted open book shown by Lisca in \cite{Lisca} (note that if we use $(S^1 \times [0,1], \tau)$ instead of $(S^1 \times [0,1], \tau^2)$ we once again get something that is not right-veering). We can see this open book in \Cref{fig:Liscaexample}.

\begin{figure}[htp]
    \centering
    \includegraphics[height=4cm]{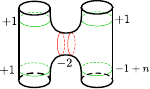}
    \caption[Overtwisted binding sums.]{Binding sum of $(S^1 \times [0,1],\tau^n)$ and $(S^1 \times [0,1], \tau^2)$.}
    \label{fig:Liscaexample}
\end{figure}

\end{example}

\section{Non-symmetric sums}
All of the examples we have studied until now involved having some boundary components of the original open books that are not used in the sum. However, even if we perform the binding sum on all binding components we might end up with a manifold that has vanishing contact invariant. For example, a manifold has Giroux torsion if and only if it can be expressed as the binding sum of three open books, two of which are $(S^1 \times [0,1], \textrm{Id})$ (see \cite{Klukas}), so if the third is a Stein fillable open book, or indeed one with non-vanishing invariant; we have a sum of manifolds with non-vanishing invariants giving a manifold with vanishing invariant (since manifolds with Giroux torsion have vanishing invariant by \cite{Ghiggini}).

It is not true, however, that the result of a binding sum is always overtwisted or has vanishing invariant. A first example is the open book on the left--hand side of \Cref{fig:S3BindingSum}. Moreover, we know by work of  Wendl and Latschev--Wendl (see \cite{Wendl, Latschev-Wendl, Klukas}) that if we sum two copies of the same open book (a \emph{symmetric sum}) with the identity as monodromy along all boundary components, the result is a Stein fillable contact structure; while summing two different surfaces (a \emph{non-symmetric sum}), still with the identity as monodromy,  along all boundary components give contact manifolds with vanishing contact invariant. 

The original result used $J$--holomorphic curves, and we want to understand this from the open book perspective; in particular by finding a chain in Heegaard Floer homology whose differential is the contact class, and using Klukas' result about the neighbourhood of the binding in the binding sum. For simplicity we work with the hat version $\hfhat$, and use $\mathbb{F}_2$ coefficients. Many of the examples are at present too complicated to handle, but we are able to give chains in an infinite family of sums, namely, summing a genus $g >0$ surface with two boundary components with an annulus, where both monodromies are the identity.

\begin{figure}[htp]
    \centering
    \includegraphics[height=6cm]{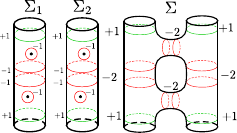}
    \caption[]{On the left, the two annuli with the navels corresponding to all boundary components. On the right, the result of the binding sum.}
    \label{fig:T^3}
\end{figure}

\begin{figure}[htp]
    \centering
    \includegraphics[height=10cm]{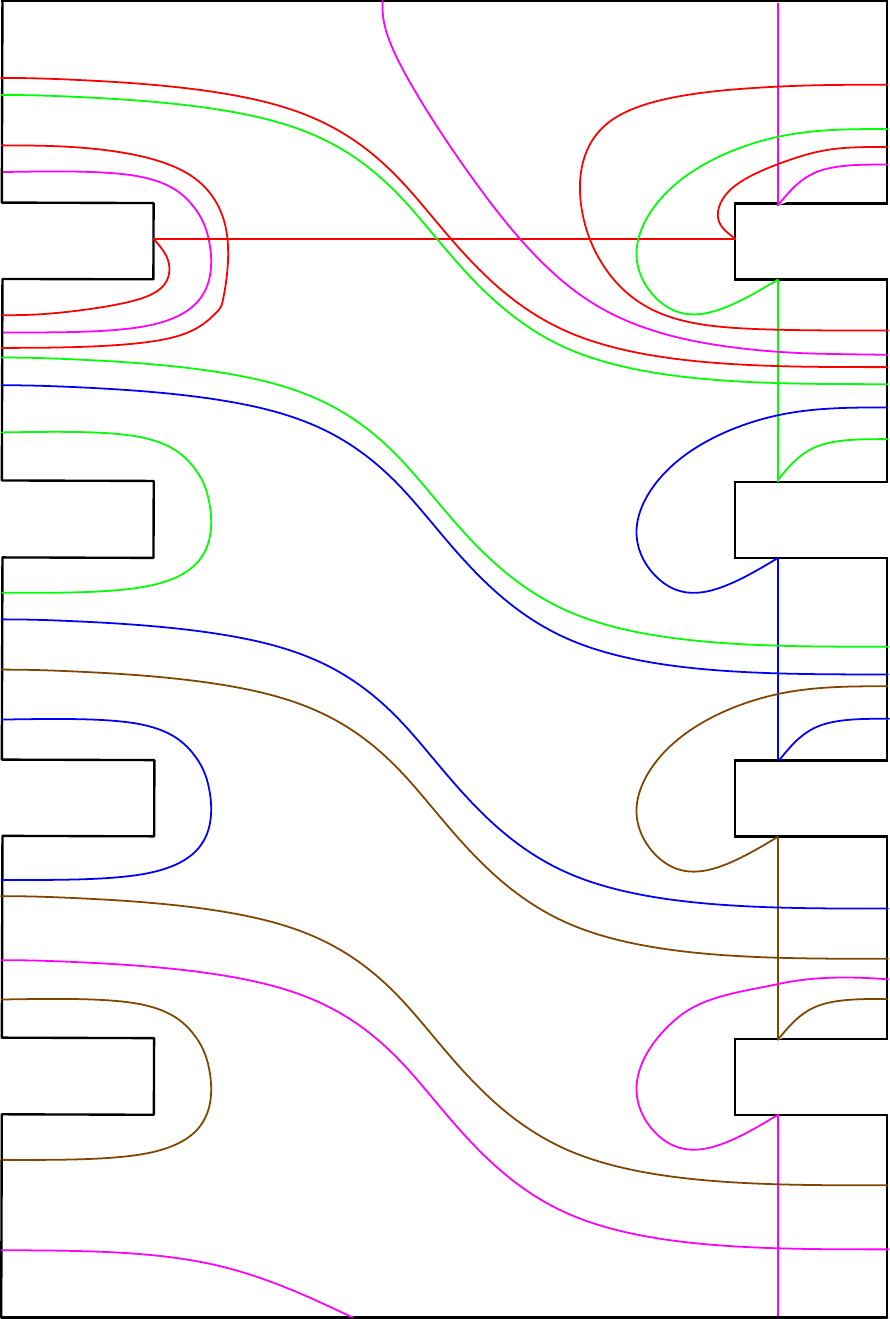}
    \caption[]{The open book for $(S^1 \times [0,1],\textrm{Id}) \boxplus (S^1 \times [0,1],\textrm{Id})$, with arcs forming a basis and their images. The top and bottom of the diagram, as well as the left and right, are identified.}
    \label{fig:T^3Monodromy}
\end{figure}

First consider the sum of two open books $(S^1 \times [0,1],\textrm{Id})$, with the sum being performed on both boundary components. The result is an open book supporting a Stein fillable contact structure on $T^3$ \cite{Klukas, VHMthesis}, and so it has non-vanishing contact invariant. The open book is given in Figure \ref{fig:T^3}, and \Cref{fig:T^3Monodromy} gives the diagram for the surface that we will use, together with a basis of arcs and its image under the monodromy.

\begin{figure}[htp]
    \centering
    \includegraphics[height=15cm]{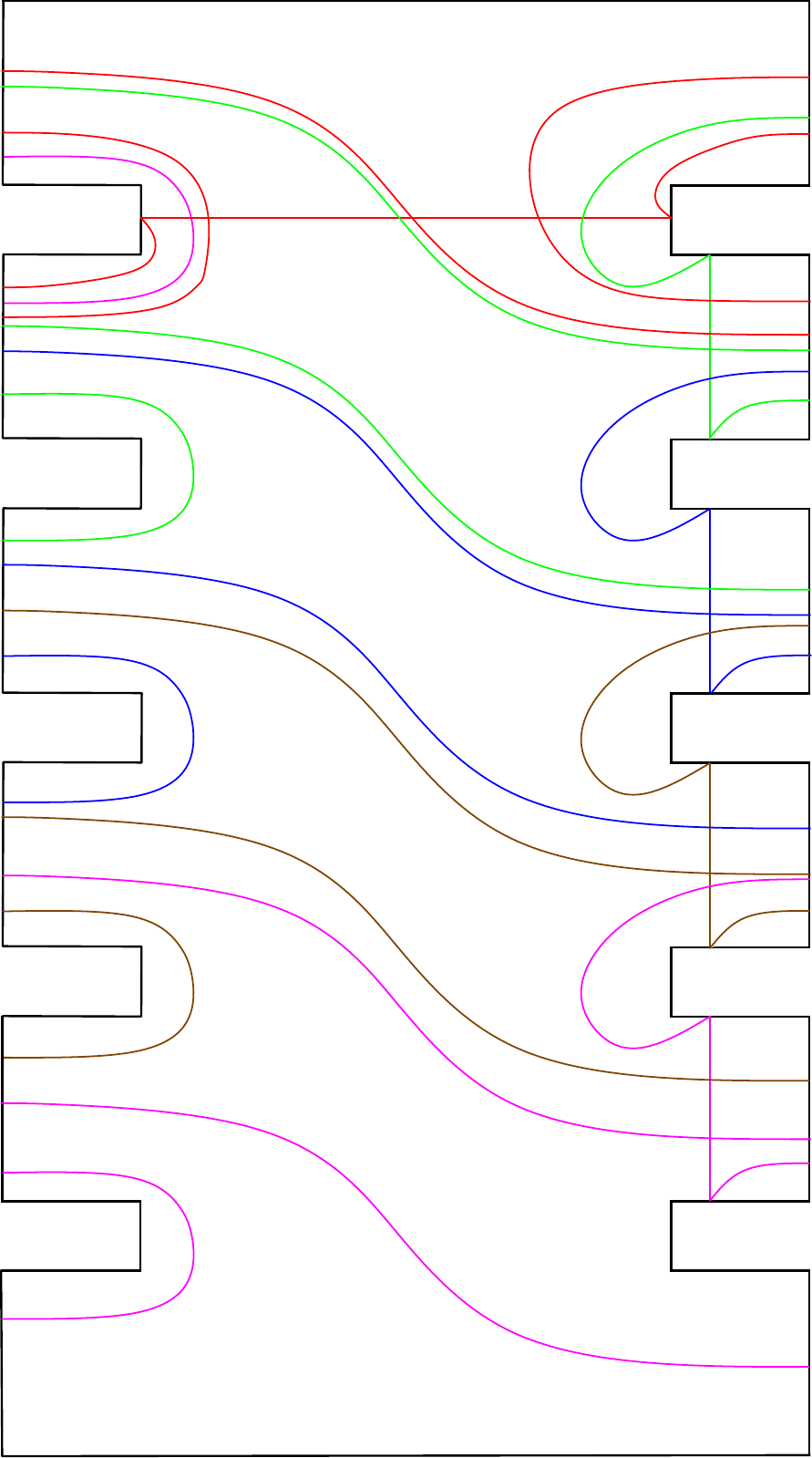}
    \caption[Partial open book for a Giroux torsion domain.]{The partial open book of a neighbourhood of a Giroux torsion domain as shown in \cite{JuhaszKang}}
    \label{fig:Giroux_torsion}
\end{figure}

Moreover, consider the partial open book with the same surface and all arcs except the one that runs through the identification between top and bottom.
Now, if we add a boundary component --more precisely, if we attach a basic slice as in \cite{Slices}--, we get the partial open book in \Cref{fig:Giroux_torsion}. This has vanishing contact class because it is the partial open book corresponding to a neighbourhood of a Giroux torsion domain (see \cite{JuhaszKang, Slices}), which, as we will calculate below, has vanishing invariant. In fact, for this specific example, we know it has vanishing invariant by \cite{Ghiggini}. 
However, the explicit calculation sidesteps technical parts about sutured Floer homology, contact surgery, and basic slices; and, more importantly, allows us to establish that the contact class vanishes in other examples, e.g. those from \Cref{thm:main}.
The next subsection presents this calculation, which is the technical heart of this paper.

\subsection{Contact class of a neighbourhood of a Giroux torsion domain}

We will adopt the convention used in \cite{JuhaszKang} for denoting generators of the Heegaard Floer complex, which is as follows. Each generator will be denoted by a tuple where the $i$-th component will be an intersection point lying in $\alpha_i$ (the arcs of the basis), where $\alpha_1$ is the horizontal arc, $\alpha_2$ is the arc immediately below $\alpha_1$ and the rest are ordered from top to bottom. The number in that component will denote the position of the intersection point in $\alpha_i$, with the order going from right to left in $\alpha_1$ and from top to bottom in the rest. In particular the contact class $c(\xi)$ is represented by the tuple $(1,1,1,1,1)$. 

\begin{figure}[htp]
    \centering
    \includegraphics[height=18cm]{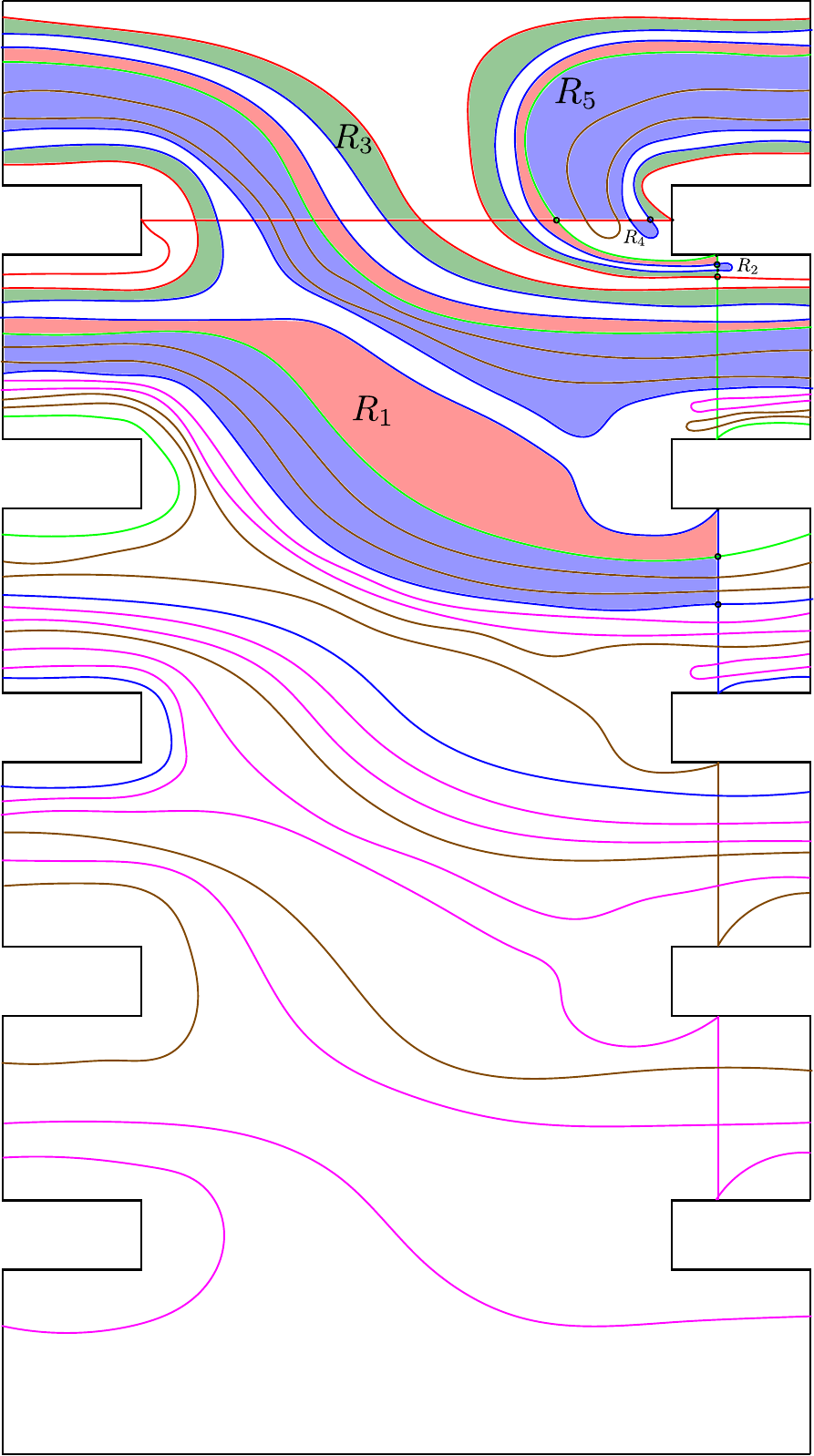}
    \caption[Nice diagram]{Nicefied diagram for the partial open book for the neighbourhood of a Giroux torsion domain, together with the first five domains.}
    \label{fig:NiceDiagram1}
\end{figure}

Now, we can apply the Sarkar-Wang algorithm \cite{NiceDiagram} to obtain a nice diagram, which is represented in \Cref{fig:NiceDiagram1}, thus turning the problem of finding the differentials into a combinatorial problem, since the only possible domains are rectangles and bigons. The chain we obtain is then the following:

There is a rectangle going from $(1,2,2,1,1)$ to $(1,1,1,1,1)$, and the only other domain coming out of this point is a bigon to $(1,3,2,1,1)$. Next, there is a rectangle from $(2,4,2,1,1)$ to $(1,3,2,1,1)$, and the only other domain coming out of it is a bigon to $(3,4,2,1,1)$. Then, there is a rectangle from $(6,4,5,1,1)$ to $(3,4,2,1,1)$, and the only other domain coming out of it is a rectangle to $(9,1,5,1,1)$. Then, there is a rectangle from $(9,1,4,2,1)$ to $(9,1,5,1,1)$, and the only other domain coming out of it is a bigon going to $(9,1,3,2,1)$. Next, there is a rectangle going from $(9,15,2,2,1)$ to $(9,1,3,2,1)$, and the only other domain coming out of it is a bigon going to $(9,14,2,2,1)$. Then there is a rectangle from $(9,11,2,5,1)$ to $(9,14,2,2,1)$, but now there are two other domains coming out of it, a rectangle going to $(16,4,2,5,1)$ and a rectangle going to $(3,5,2,5,1)$. Now there is a bigon from $(2,5,2,5,1)$ that goes to $(3,5,2,5,1)$, and the only other domain coming out of it is a rectangle going to $(1,6,2,5,1)$. We then consider the point $(17,4,2,5,1)$, which has a bigon to $(16,4,2,5,1)$ and a rectangle to $(1,6,2,5,1)$, but also a rectangle to $(18,3,2,5,1)$. Next, there is a bigon from $(18,2,2,5,1)$ to $(18,3,2,5,1)$, and the only other domain coming out of it is a rectangle to $(18,1,1,5,1)$. Then, there is a rectangle from $(18,1,1,4,2)$ to $(18,1,1,5,1)$, and the only other domain coming out of it is a bigon to $(18,1,1,3,2)$. Next, there is a rectangle from $(18,1,10,2,2)$ to $(18,1,1,3,2)$ and the only other domain coming out of it is a bigon to $(18,1,9,2,2)$. Finally, there is a rectangle from $(18,1,5,6,2)$ to $(18,1,9,2,2)$ and there are no other domains coming out of this point. We gather the domains in \Cref{table:Jplus} for clarity. Thus, by working with $\mathbb{F}_2$ coefficients we get that

$ \partial( (1,2,2,1,1) + (2,4,2,1,1) + (6,4,5,1,1) + (9,1,4,2,1) + (9,15,2,2,1) +  (9,11,2,5,1) + (2,5,2,5,1) + (17,4,2,5,1) + (18,2,2,5,1) + (18,1,1,4,2) + (18,1,10,2,2) + (18,1,5,6,2)) = (1,1,1,1,1)$

To confirm this, observe that every point except the contact generator appears twice in the differential--checking can be done in \Cref{table:Jplus}.

\begin{figure}[htp]
    \centering
    \includegraphics[height=18cm]{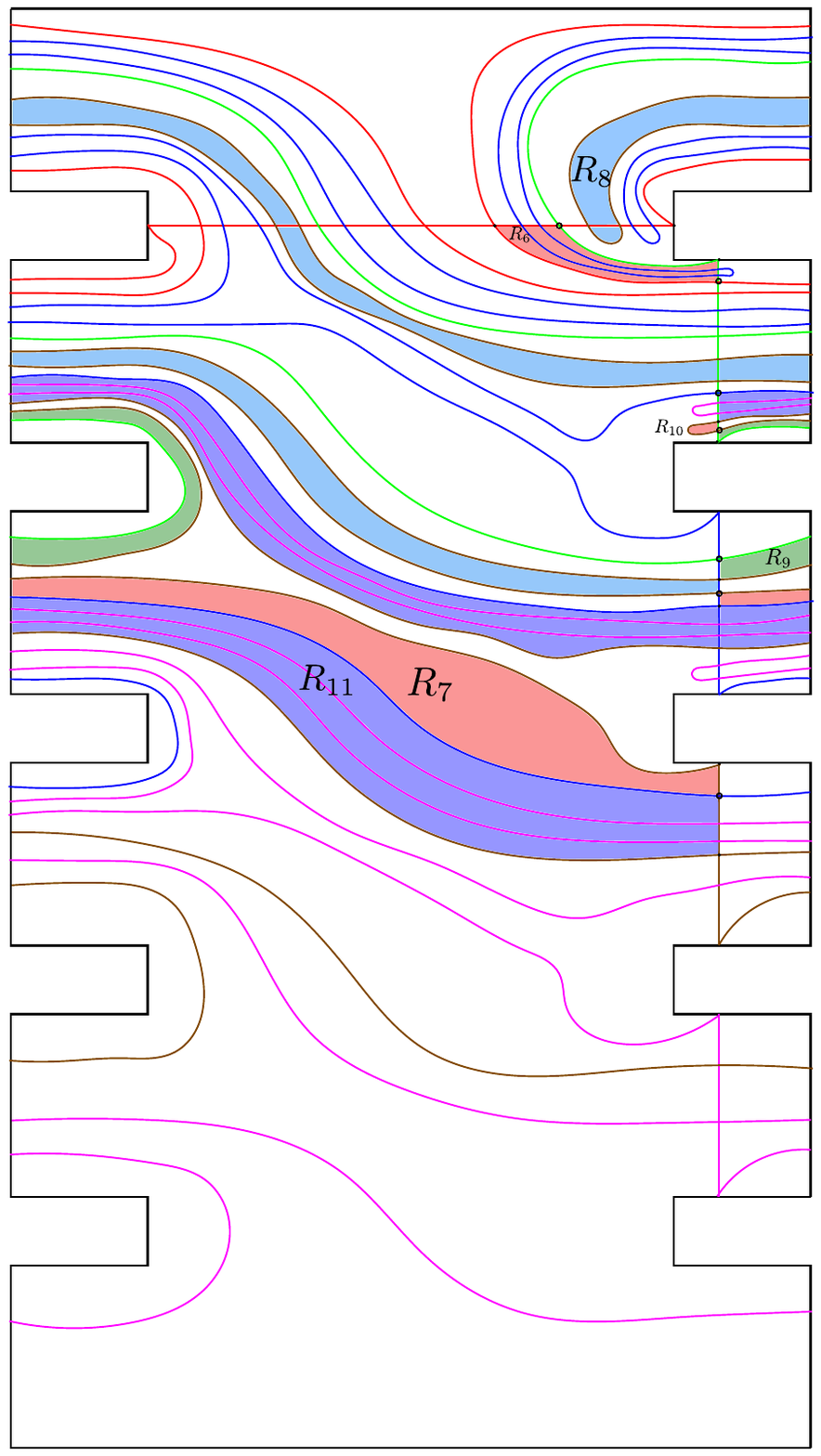}
    \caption[Nice diagram]{The domains $R_6$ to $R_{11}$.}
    \label{fig:NiceDiagram2}
\end{figure}

\begin{figure}[htp]
    \centering
    \includegraphics[height=18cm]{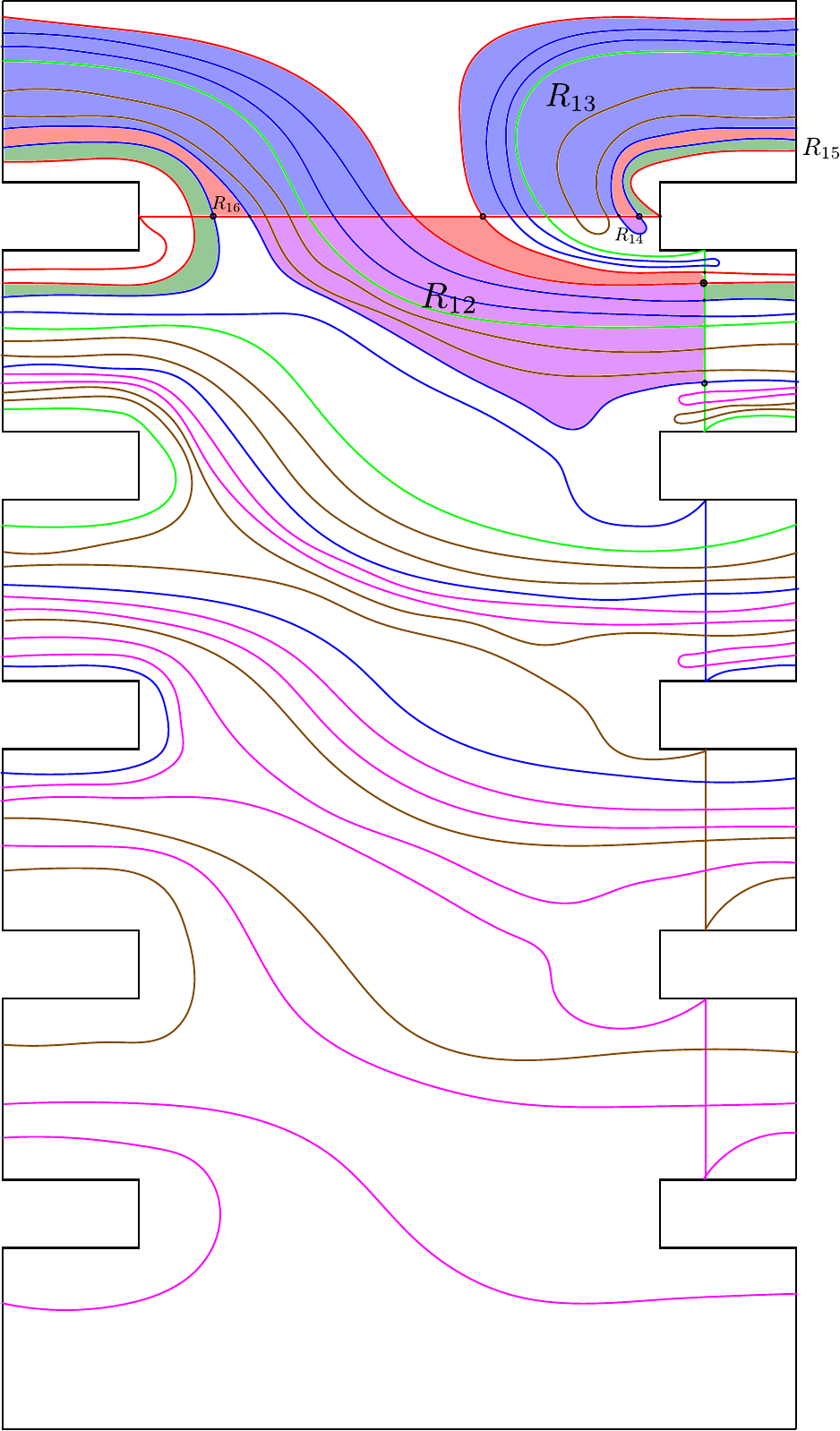}
    \caption[Nice diagram]{The domains $R_{12}$ to $R_{16}$. Note that domains $R_{12}$ (red) and $R_{13}$ (blue) overlap, this is illustrated by shading in purple. Moreover, the bigon $R_{14}$ is contained in the bigon $R_{16}$ (red), this is also illustrated by shading it in purple. }
    \label{fig:NiceDiagram3}
\end{figure}

\begin{figure}[htp]
    \centering
    \includegraphics[height=18cm]{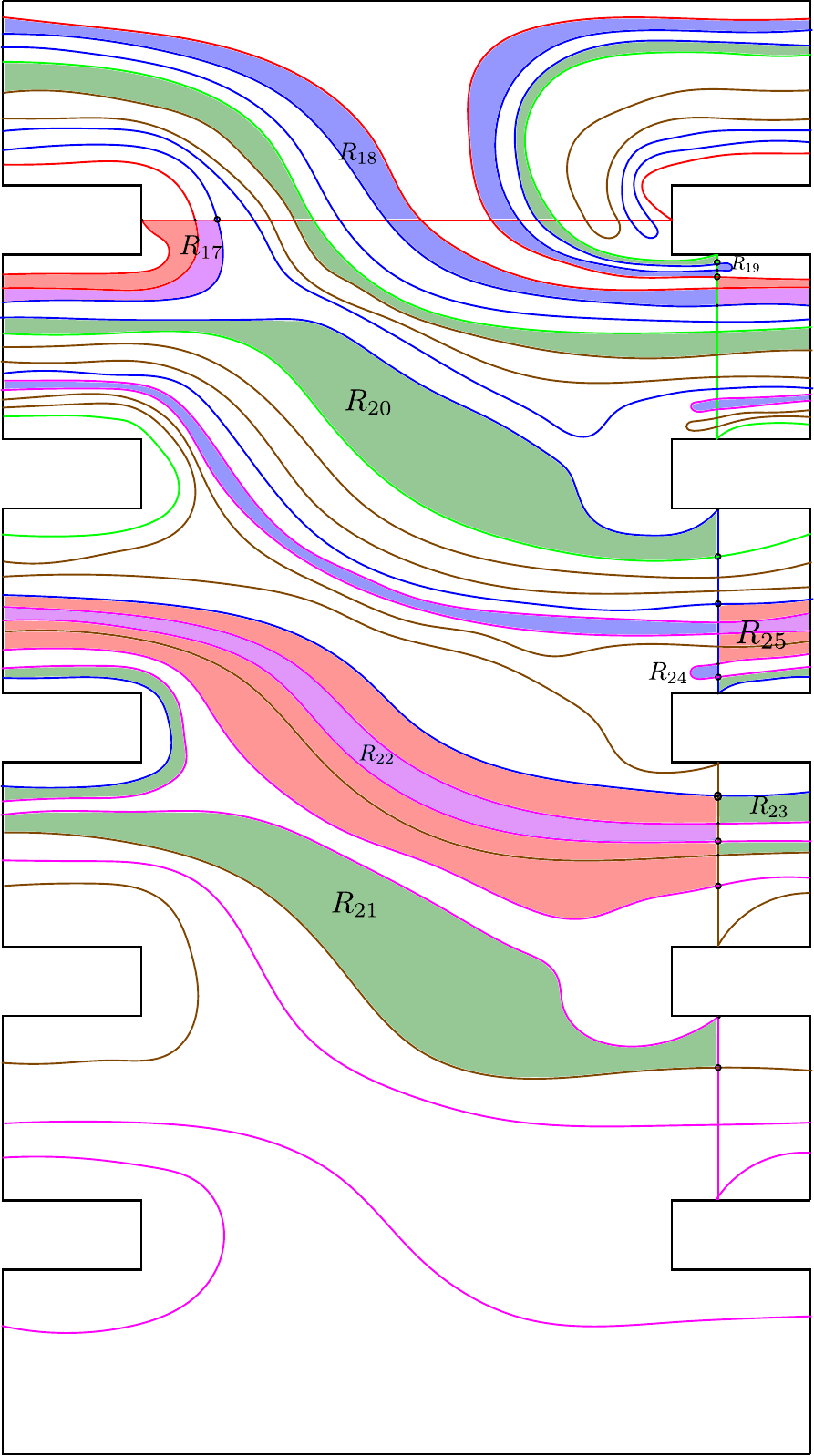}
    \caption[Nice diagram]{The domains $R_{17}$ to $R_{25}$. As before, when blue and red domains overlap, this is illustrated by shading the overlap in purple.}
    \label{fig:NiceDiagram4}
\end{figure}

We can see the domains in this chain in \Cref{fig:NiceDiagram1}, \Cref{fig:NiceDiagram2}, \Cref{fig:NiceDiagram3}, and \Cref{fig:NiceDiagram4}. The intersection points that are the sources of a domain are labeled as $\circ$-points and the points that are the targets are labeled as $\bullet$-points. Some of the domains on the surface are repeated. However, since the other intersection points in the tuple are different, they actually correspond to distinct $J$--holomorphic disks.

\begin{remark} \label{rmk:computerassistance}
    This computation was done assisted by the computer. While it is easy to check that every domain given here is indeed a domain involved in the differential, due to the large number of intersection points it is not easy to determine that there are no more domains. Moreover, the chain we obtain is different from the one computed in \cite{JuhaszKang}. However, as the aim of this section is to prove \Cref{thm:main}, we leave this discussion for \Cref{sec:spectralorder}.
\end{remark}

\begin{table} [ht]  
\centering
\begin{tabular}{ |c|c|c|c| }
 \hline
 From & To & $J_+$ & Region name \\ 
 \hline
 $(1,2,2,1,1)$ & $(1,1,1,1,1)$ & 0 & $R_1$\\ 
 $(1,2,2,1,1)$ & $(1,3,2,1,1)$ & 0 &  $R_2$\\ 
 $(2,4,2,1,1)$ & $(1,3,2,1,1)$ & 0 & $R_3$\\
 $(2,4,2,1,1)$ & $(3,4,2,1,1)$ & 0 & $R_4$\\
 $(6,4,5,1,1)$ & $(3,4,2,1,1)$ & 2 & $R_5$\\
 $(6,4,5,1,1)$ & $(9,1,5,1,1)$ & 0 & $R_6$\\
 $(9,1,4,2,1)$ & $(9,1,5,1,1)$ & 0 & $R_7$\\
 $(9,1,4,2,1)$ & $(9,1,3,2,1)$ & 0 & $R_8$\\
 $(9,15,2,2,1)$ & $(9,1,3,2,1)$ & 0 & $R_9$\\
 $(9,15,2,2,1)$ & $(9,14,2,2,1)$ & 0 & $R_{10}$\\
 $(9,11,2,5,1)$ & $(9,14,2,2,1)$ & 2 & $R_{11}$\\
 $(9,11,2,5,1)$ & $(16,4,2,5,1)$ & 2 & $R_{12}$\\
 $(9,11,2,5,1)$ & $(3,5,2,5,1)$ & 2 & $R_{13}$\\
 $(2,5,2,5,1)$ & $(3,5,2,5,1)$ & 0 & $R_{14}$\\
 $(2,5,2,5,1)$ & $(1,6,2,5,1)$ & 0 & $R_{15}$\\
 $(17,4,2,5,1)$ & $(16,4,2,5,1)$ & 0 & $R_{16}$\\
 $(17,4,2,5,1)$ & $(1,6,2,5,1)$ & 0 & $R_{17}$\\
 $(17,4,2,5,1)$ & $(18,3,2,5,1)$ & 0 & $R_{18}$\\
 $(18,2,2,5,1)$ & $(18,3,2,5,1)$ & 0 & $R_{19}$\\
 $(18,2,2,5,1)$ & $(18,1,1,5,1)$ & 0 & $R_{20}$\\
 $(18,1,1,4,2)$ & $(18,1,1,5,1)$ & 0 & $R_{21}$\\
 $(18,1,1,4,2)$ & $(18,1,1,3,2)$ & 0 & $R_{22}$\\
 $(18,1,10,2,2)$ & $(18,1,1,3,2)$ & 0 & $R_{23}$\\
 $(18,1,10,2,2)$ & $(18,1,9,2,2)$ & 0 & $R_{24}$\\
 $(18,1,5,6,2)$ & $(18,1,9,2,2)$ & 2 & $R_{25}$\\
 \hline
\end{tabular}
\caption{The domains in the chain killing the contact class together with their $J_+$, which will be necessary for the spectral order computation in \Cref{sec:spectralorder}. The region names correspond to the labels in the figures.}
\label{table:Jplus}
\end{table}

\begin{figure}[ht]
    \centering
    \includegraphics[height=13cm]{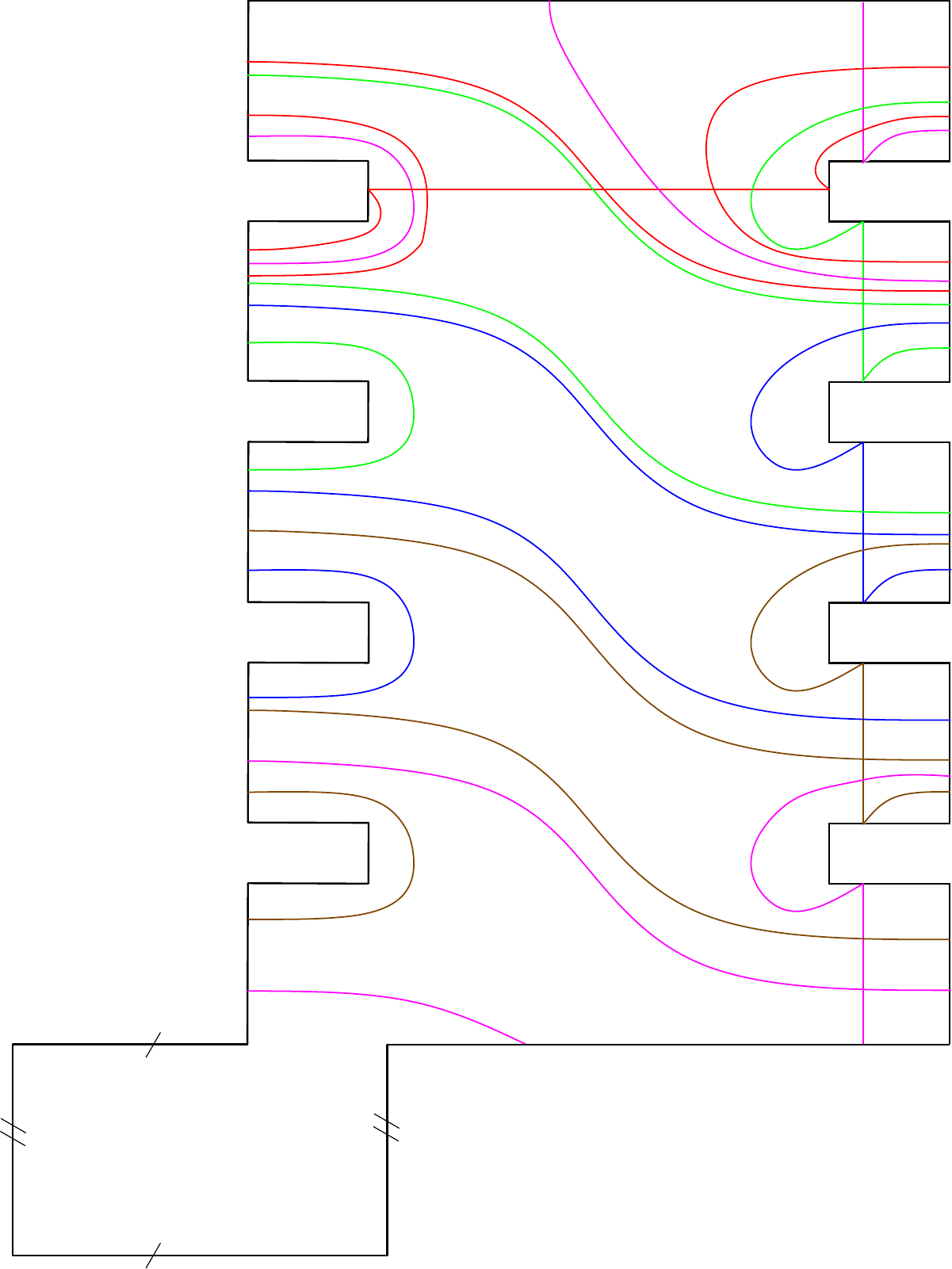}
    \caption[Open book of a non-symmetric binding sum.]{An open book for the non-symmetric sum $(\Sigma_{1,2},\textrm{Id}) \boxplus (A,\textrm{Id})$, where the identifications are as before but we also need to identify the boundary of the rectangle at the bottom left of the figure, which we do in a standard way for a torus by identifying top and bottom and left and right. The two extra arcs needed to complete to a basis are omitted for simplicity, since they are not used in the computation.}
    \label{fig:Non_symmetric}
\end{figure}

We remark that we can see the chain killing the contact class in the original diagram, after undoing the finger moves that created the nice diagram. Using the same notational convention for \Cref{fig:Giroux_torsion} as we did before, it is the following:

There is an embedded annulus from (2,2,3,1,1) to (1,1,1,1,1). The only other domain coming out of (2,2,3,1,1) is a disk going to (3,1,3,1,1). Next, there is an embedded annulus from (3,6,2,3,1). The only other domain coming out of (3,6,2,3,1) is an immersed annulus going to (5,1,1,3,1). Finally, there is an embedded annulus going from (5,1,3,4,2) to (5,1,1,3,1) and no other domains come out of (5,1,3,4,2). This means that 

$ \partial ((2,2,3,1,1) + (3,6,2,3,1) + (5,1,3,4,2)) = (1,1,1,1,1) + 2(3,1,3,1,1) + 2(5,1,1,3,1) = (1,1,1,1,1) $  because we are working with $\mathbb{F}_2$ coefficients.

\subsection{Proof of \Cref{thm:main}}

Using the chain computed above, we can show an explicit chain that causes vanishing of the contact class for an infinite family of binding sums. Thus, we are ready to prove \Cref{thm:main}.

\begin{proof} [Proof of \Cref{thm:main}]
We begin by proving part $a)$. Observe that the domains in this partial diagram will also be in the diagram for $T^3$, since we can see the last boundary component as the first one by an identification. The reason why the contact invariant does not vanish in this open book is that there are more domains coming from the rest of the arcs in the basis once we have gone back to the top. However, these are the only extra domains. This means that, if we do not use domains that require the identification of top and bottom, we can get a chain killing the contact class in the partial case. Now, adding genus on the bottom part of the open book means that there will be no domains that use both top and bottom. Thus, the chain we used before also kills the contact invariant in this case. However, now observe that adding genus to the open book in this manner amounts to adding genus to one of the summands in the binding sum, and we can add arbitrarily large genus. \Cref{fig:Non_symmetric} shows the case where the genus is $1$, and we can see that all the domains that we computed before are present in this open book, and there are no new ones. This proves $a)$.

For part $b)$, observe that a diffeomorphism of $\Sigma_g$ does not change the images of the arcs involved in the computation of the chain from before, see for example \Cref{fig:Non_symmetric}. Therefore, it still kills the contact class.

\end{proof}

 \label{sec:nonsym}

\section{Spectral order of a neighbourhood of a Giroux torsion domain} \label{sec:spectralorder}
Up until now most of the work presented was done in the author's thesis \cite{mythesis}. In doing the computation for the contact class in the proof of \Cref{thm:main}, we noticed that \cite{JuhaszKang} made a similar calculation. However, the two chains differ, and we can show that there is a subtle mistake in \cite{JuhaszKang}. Indeed, there the chain ends with the point $(9,11,2,3,2)$, with the only domain coming out of it being a rectangle to $(9,12,2,2,2)$. However, there are also rectangles going to $(16,4,2,3,2)$ and $(3,5,2,3,2)$, and so this chain does not work. As pointed out in \Cref{rmk:computerassistance}, it is easy to check whether a given domain contributes to the differential but not as easy to find all domains, which is why computer assistance was useful for this computation.

We remark that it was shown in \cite{Ghiggini} that the contact class vanishes in this case, so this does not provide a new result. However, Juh\'{a}sz and Kang use the computation of vanishing of the contact class to provide an upper bound for the spectral order, which is a refinement of the contact class defined in \cite{SpectralOrder}. In particular, they obtain \Cref{cor:spectralbound}. We refer to \cite{SpectralOrder, JuhaszKang} for definition of the invariant and how to compute upper bounds for it, and simply note that the main input we need is that $\partial_i$ counts $J$--holomorphic disks with $J_+ = 2i$.

With our corrected chain, we can prove \Cref{cor:spectralbound}.  

\begin{proof} [Proof of \Cref{cor:spectralbound}]
    
Recall from before that the chain killing the contact class is 

$ \partial( (1,2,2,1,1) + (2,4,2,1,1) + (6,4,5,1,1) + (9,1,4,2,1) + (9,15,2,2,1) +  (9,11,2,5,1) + (2,5,2,5,1) + (17,4,2,5,1) + (18,2,2,5,1) + (18,1,1,4,2) + (18,1,10,2,2) + (18,1,5,6,2)) = (1,1,1,1,1)$
\medskip

Now let $b_0 = (1,2,2,1,1) + (2,4,2,1,1) + (6,4,5,1,1) + (9,1,4,2,1) + (9,15,2,2,1) + (2,5,2,5,1) + (17,4,2,5,1) + (18,2, 2, 5, 1) + (18,1,1,4,2) + (18,1,10,2,2)$,

$b_1 = (6,4,5,1,1) + (9,1,4,2,1) + (9,15,2,2,1) + (9,11,2,5,1) + (18,1,5,6,2) + (2,5,2,5,1) + (17,4,2,5,1) + (18,2,2,5,1) + (18,1,1,4,2) + (18,1,10,2,2) + (18,1,5,6,2)$, and

$b_2 = (9,11,2,5,1) + (18,1,5,6,2)$.

Then, from \Cref{table:Jplus} we get the following differentials (recall we work over $\mathbb{F}_2$).

$\partial_0 b_0 = (1,1,1,1,1) + (9,14,2,2,1) + (16,4,2,5,1) + (3,5,2,5,1) + (18,1,9,2,2)$,

$\partial_0 b_1 = (9,14,2,2,1) + (16,4,2,5,1) + (3,5,2,5,1) + (18,1,9,2,2)$ and $\partial_1 b_1 = (9,14,2,2,1) + (16,4,2,5,1) + (3,5,2,5,1) + (18,1,9,2,2)$,

$\partial_0 b_2 = \partial_2 b_2 = 0$ and $\partial_1 b_2 = (9,14,2,2,1) + (16,4,2,5,1) + (3,5,2,5,1) + (18,1,9,2,2) = \partial_0 b_1$.

Therefore we get that $ \partial_0 b_0 + \partial_1 b_1 + \partial_2 b_2 = (1,1,1,1,1)$, $\partial_0 b_1 + \partial_1 b_2 = 2\partial_0 b_1 = 0$, and $\partial_0 b_2 = 0$.
Thus, setting $b_i = 0 $ for $i > 2$, $c_i = 0$ for  $i\geq1$, and $c_0 = (1,1,1,1,1)$, we obtain that $\hat{\partial} (b_i)_{i \in \mathbb{N}} = (c_i)_{i \in \mathbb{N}}$. Following \cite{SpectralOrder,JuhaszKang}, this is enough to guarantee the result.

\end{proof}

It is conjectured--see \cite{SpectralOrder} for the closed case, \cite{JuhaszKang} for the case of manifolds with convex boundary--that the spectral order of a manifold with Giroux torsion is $1$. However, the upper bound of $2$ cannot be improved using the chain computed in this paper. Indeed, to obtain a bound of $1$, we would need $b_0$ and $b_1$ such that $\partial_0 b_0 + \partial_1 b_1 = (1,1,1,1,1)$ and $\partial_0 b_1 = 0$. But then $(6,4,5,1,1)$--with nontrivial $\partial_0$--must appear in $b_1$, since it also has nontrivial $\partial_1$. But then by the condition $\partial_0 b_1 = 0$ we also need $(9,1,4,2,1)$ in $b_1$, which in turn forces $(9,15,2,2,1)$ to be in $b_1$. However, $\partial_0 (9,15,2,2,1)$ contains $(9,14,2,2,1)$, and there are no domains going to it with $J_+ = 0$, so $\partial_0 b_1$ cannot vanish.

\bibliographystyle{myamsalpha} 
\bibliography{References}

\end{document}